# Equivalence of time-domain inverse problems and boundary spectral problems

A. Katchalov [*]   Y. Kurylev[†]   M. Lassas[‡]   N. Mandache[§]

October 31, 2018

**Abstract.** We consider inverse problems for wave, heat and Schrödinger-type operators and corresponding spectral problems on domains of $\mathbf{R}^n$ and compact manifolds. Also, we study inverse problems where coefficients of partial differential operator have to be found when one knows how much energy it is required to force the solution to have given boundary values, i.e., one knows how much energy is needed to make given measurements. The main result of the paper is to show that all these problems are shown to be equivalent.

## 1 Introduction

The goal of this paper is to explain the equivalence of various types of data used in inverse boundary value problems. Our main interest lies in the inverse boundary spectral problems and the inverse boundary value problems in the time-domain for equations with time-independent coefficients. There are numerous examples when a solution of a particular inverse problem is used to solve an inverse problem of another type. For example, A. Nachman, J. Sylvester, and G. Uhlmann [NSU] solved the inverse boundary spectral problem for a Schrödinger operator by reducing it to the inverse boundary value problem in fixed frequency and then using the method of complex geometric optics [SU]. Similarly, inverse problems in the time-domain are often

[*]Steklov Mathematical Institute, RAN, Fontanka 27, 191011, St. Petersburg, Russia
[†]Department of Mathematical Sciences, Loughborough University, Loughborough, LE11 3TU, UK
[‡]Rolf Nevanlinna Institute, University of Helsinki Helsinki, P.O.Box 4, FIN-00014, Finland
[§]Department of Mathematical Sciences, Loughborough University, Loughborough, LE11 3TU, UK

reduced to problems in the frequency domain (see e.g. [Is]). In this paper we are interested in the equivalence of different types of boundary data for inverse problems and describe some procedures for transforming between these problems. We believe that the equivalence of various inverse problems is useful from the theoretical point of view since it enables us to translate the results obtained for one inverse problem to other problems. The equivalence is also useful for applications as it makes possible to use reconstruction algorithms developed for some particular type of inverse problems to different inverse problems. Moreover, in various applications some measurements are considerably more difficult to make when compared to other measurements which, however, can give an equivalent information. For instance, it is often difficult to measure the phase of a physical field but the energy which is required to force the boundary value of the field to a given one is often known very precisely. In theoretical inverse problems this idea goes back to A. Calderón who in his seminal paper of 1980 considered an inverse problem for the conductivity equation from the point view of energy measurements. Similarly, the energy or interference based measurements have been used in many applications, e.g. in impedance tomography (see e.g. [CIN]) and near field optical tomography (see e.g. [SM]).

We start with an elliptic problem in the frequency domain. Let $\Omega \subset \mathbf{R}^m$ be a smooth bounded domain and $a(x, D)$ be an elliptic 2nd-order partial differential operator,

$$a(x,D)v(x) = -g^{-1/2}\partial_j\left(g^{1/2}g^{jk}(x)\partial_k v(x)\right) + q(x)v(x), \tag{1}$$

where $[g^{jk}(x)]$ is a positive definite smooth real matrix in $\overline{\Omega}$, $[g_{ij}(x)]$ is the inverse matrix of $[g^{jk}(x)]$, $g(x) = \det(g_{ij})$, and $q(x)$ is a smooth real function. As usual we use Einstein's summation over repeated upper and lower indices. The operator $a(x, D)$ is a Schrödinger operator in $\Omega$ which corresponds to the metric tensor $g_{ij}$.

Two objects related to the Dirichlet problem are the boundary spectral data and the Dirichlet-to-Neumann map $\Lambda_z$. Namely, if $A$ is the operator

$$Av = a(x,D)v$$

with

$$\mathcal{D}(A) = H^2(\Omega) \cap H_0^1(\Omega), \tag{2}$$

then its boundary spectral data is the collection

$$\{\lambda_l,\ B\varphi_l|_{\partial\Omega}:\ \ l = 1, 2, \dots\}, \tag{3}$$

where $\lambda_l$ and $\varphi_l$ are the eigenvalues and the normalized eigenfunctions of the operator $A$,

$$a(x, D)\varphi_l = \lambda_l \varphi_l, \quad \varphi_l|_{\partial\Omega} = 0.$$

The boundary operator $B$ is given by

$$B\varphi_l = g^{1/2} g^{jk} \nu_j \partial_k \varphi_l, \tag{4}$$

where $\nu = (\nu_1, \ldots, \nu_m)$ is the interior unit normal vector to $\partial\Omega$ in the Euclidean metric. To clarify the meaning of the boundary operator $B$, we note that if $n$ is the interior unit normal vector to $\partial\Omega$ with respect to the metric $g$, then $Bu = \rho \partial_n u$, where the weight $\rho$ is given by

$$\rho(x) = \det(g_{\partial\Omega}(x))^{1/2} (g^{ij}(x)\nu_i \nu_j)^{1/2}. \tag{5}$$

Here $g_{\partial\Omega}$ is the restriction of the metric tensor $g_{jk}(x)$ to the tangent plane $H = T_x(\partial\Omega)$ of $\partial\Omega$ at $x$. Thus, integration by parts gives rise to the formula

$$\int_\Omega a(x,D) v\, \overline{u}\, dV_g = \int_\Omega \left( g^{ij} \partial_j v\, \overline{\partial_k u} + qv\overline{u} \right) dV_g + \int_{\partial\Omega} Bv|_{\partial\Omega}\, \overline{u|_{\partial\Omega}} dS.$$

Here $dS$ is the Euclidean area element on $\partial\Omega$ while $dV_g = g^{1/2} dx$ is the Riemannian volume element on $\Omega$ in the metric $g$.

Then, for $z \neq \lambda_l$, the (fixed-frequency) Dirichlet-to-Neumann map $\Lambda_z$ is defined by

$$\Lambda_z h = B v_z^h|_{\partial\Omega}, \quad h \in C^\infty(\partial\Omega), \tag{6}$$

where $v_z^h$ solves the boundary value problem

$$a(x, D) v_z^h = z v_z^h, \quad v_z^h|_{\partial\Omega} = h. \tag{7}$$

We are interested in the case when $\Lambda_z$ is given at all possible values of the spectral parameter $z \in \mathbf{C}$. This data is called the Gel'fand spectral data.

Operator $a(x, D)$ defines in a natural way initial boundary value problems for the wave equation,

$$\begin{aligned}(\partial_t^2 + a(x,D))u^f &= 0 \quad \text{in} \quad \Omega \times \mathbf{R}_+, \\ u^f|_{\partial\Omega \times \mathbf{R}_+} = f \in C_0^\infty(\partial\Omega \times \mathbf{R}_+), \quad & u^f|_{t=0} = \partial_t u^f|_{t=0} = 0;\end{aligned} \tag{8}$$

for the heat equation,

$$\begin{aligned}(\partial_t + a(x,D))w^f &= 0 \quad \text{in} \quad \Omega \times \mathbf{R}_+, \\ w^f|_{\partial\Omega \times \mathbf{R}_+} &= f, \quad w^f|_{t=0} = 0;\end{aligned} \tag{9}$$

and for the non-stationary Schrödinger equation,

$$(i\partial_t + a(x, D))\psi^f = 0 \quad \text{in } \Omega \times \mathbf{R}_+, \tag{10}$$
$$\psi^f|_{\partial\Omega \times \mathbf{R}_+} = f, \quad \psi^f|_{t=0} = 0.$$

For the initial boundary value problems (8)-(10) we define the non-stationary Dirichlet-to-Neumann maps (response operators) $R^w$, $R^h$, and $R^s$ which correspond to the wave, heat, and non-stationary Schrödinger equations,

$$R^w f = Bu^f|_{\partial\Omega \times \mathbf{R}_+}, \tag{11}$$
$$R^h f = Bw^f|_{\partial\Omega \times \mathbf{R}_+}, \tag{12}$$
$$R^s f = B\psi^f|_{\partial\Omega \times \mathbf{R}_+}. \tag{13}$$

We note that $R^w$ and $R^h$ are real-valued operators, i.e. $R^w f$, $R^h f$ are real as soon as $f$ is real. Later, this will allow us to restrict ourselves to the consideration of only real $f's$ when dealing with problems (8) and (9).

There is a natural concept of energy for the wave equation given by

$$E^w(u, t) = \frac{1}{2} \int_\Omega \left( g^{jk} \partial_j u(x, t) \partial_k \overline{u(x, t)} + q(x)|u(x, t)|^2 + |\partial_t u(x, t)|^2 \right) dV_g(x).$$

Similarly, for the non-stationary Schrödinger equation we define the energy

$$E^s(\psi, t) = \frac{1}{2} \int_\Omega \left( g^{jk} \partial_j \psi(x, t) \partial_k \overline{\psi(x, t)} + q(x)|\psi(x, t)|^2 \right) dV_g(x).$$

The notion of energy (or, better say, heat content) for the heat equation has sense in the case $q = 0$ and is given by

$$E^h(w, t) = \int_\Omega w(x, t) \, dV_g(x), \tag{14}$$

with the term $qw$ physically corresponding to a heat source. For the wave and non-stationary Schrödinger equations the energy is preserved in the absence of the boundary and interior sources. In the case of the initial boundary value problems with $f \in C_0^\infty(\partial\Omega \times [0, T])$ this means that the energies $E^w(u^f, t)$ and $E^s(\psi^f, t)$ are constant for $t > T$. Therefore, the energies $E^w(u^f, t)$ and $E^s(\psi^f, t)$ are brought into $\Omega$ through the boundary $\partial\Omega \times \mathbf{R}_+$. We define the total energy fluxes through the boundary, $\Pi^w(f)$ and $\Pi^s(f)$ for the wave and non-stationary Schrödinger equations as

$$\Pi^w(f) = \lim_{t \to \infty} E^w(u^f, t), \quad \Pi^s(f) = \lim_{t \to \infty} E^s(\psi^f, t), \tag{15}$$

where $f \in C_0^\infty(\partial\Omega \times \mathbf{R}_+)$. Actually, $\Pi^w(f)$ and $\Pi^s(f)$ are given by some quadratic forms of $f$ (see formulae (28) and (35)). The total energy fluxes are important concepts, as they are often easier to measure in practice than the corresponding non-stationary Dirichlet-to-Neumann maps.

Now we are in the position to formulate various inverse problems in the domain $\Omega$ which are related to the above concepts.

**Elliptic problems:**

i. Given boundary spectral data, $\{\lambda_l,\ B\varphi_l|_{\partial\Omega}:\ l=1,2,\dots\}$, determine $g^{jk}$ and $q$.

ii. Given Gel'fand spectral data, $\{\Lambda_z,\ z \in \mathbf{C}\}$, determine $g^{jk}$ and $q$.

**Hyperbolic problems:**

iii. Given a hyperbolic Dirichlet-to-Neumann map $R^w$, determine $g^{jk}$ and $q$.

iv. Given a hyperbolic energy flux $\Pi^w$, determine $g^{jk}$ and $q$.

**Parabolic problems:**

v. Given a parabolic Dirichlet-to-Neumann map $R^h$, determine $g^{jk}$ and $q$.

**Non-stationary Schrödinger problems:**

vi. Given a Dirichlet-to-Neumann map for a non-stationary Schrödinger equation, $R^s$ determine $g^{jk}$ and $q$.

vii. Given an energy flux $\Pi^s$ for a Schrödinger equation, determine $g^{jk}$ and $q$.

Our main result is:

**Theorem 1** *Inverse problems i.-vi. are equivalent, i.e. any of the data i.-vi. determine all other data. In the case of problem vii. we should require, in addition, that $\lambda = 0$ is not an eigenvalue of $A$. Then vii. is equivalent to $i. - vi.$*

We note that if $0 \in \sigma(A)$, where $\sigma(A) = \{\lambda_1 < \lambda_2 \leq \cdots\}$ is the spectrum of $A$, then $\Pi^s$ determines all $\{\lambda_l, B\varphi_l|_{\partial\Omega}\}$ for $\lambda_l \neq 0$. Similarly, $\Pi^s$ in this case determines all Dirichlet-to-Neumann maps (6), (11)-(13) upto a finite-dimensional operator.

The authors want to emphasize that some of the above equivalences are known at least on the formal level. Rigorous proofs of some of them are given in [KKL], [Is]. In these cases, we will describe the corresponding results very briefly and refer to [KKL]. However, in particular for the Dirichlet problem considered in this paper, care should be taken to make these formal constructions rigorous.

Although solving inverse problems $i.$-$vii.$ is not the main goal of this paper, we would like to mention the following application of Theorem 1:

**Theorem 2** *Assume that $\Omega$ and any of the boundary data $i.$-$vii.$ are given. Then it is possible to construct $g^{jk}$ and $q$ upto a diffeomorphism, that is, we can find the equivalence class*

$$[a(x, D)] = \{\Phi^*(a(x, D)) \; : \; \Phi : \Omega \to \Omega \text{ is a diffeomorphism}, \; \Phi|_{\partial\Omega} = id|_{\partial\Omega}\},$$

*where $\Phi^*$ is the pull-back of the diffeomorphism $\Phi$.*

For problem $i.$ the proof of Theorem 2 can be found in [K1], [K2], and, in more detail, in [KKL], Section 4.5. Moreover, it is shown in [KK] that Theorem 2 is valid for the case when we know all but a finite number of boundary spectral data. Other cases follow then from Theorem 1 and the remark after it.

In Sections 2-4 we will describe some constructive methods to transform data $i.$-$vii.$ into each other. These methods require neither solving the corresponding inverse problem nor using an analytic continuation. We will not provide all technical details by either referring to [KKL] or Appendix 1. In Section 5 we will discuss some generalizations of problems $i.$-$vii.$. In particular, we consider the case when the data is given on a finite time interval or $a(x, D)$ is different from (1) with some technical details concerning the energy flux for the wave equation considered in Appendix 2.

## 2 Equivalence of the boundary spectral and Gel'fand data.

$\underline{i \to ii}$. Given $\{\lambda_l, \; B\varphi_l|_{\partial\Omega}, \; l = 1, 2, \dots\}$ one can formally construct the Dirichlet-to-Neumann map $\Lambda_z$,

$$(\Lambda_z h)(x) \underset{formally}{=} \sum_{l=1}^{\infty} \frac{B\varphi_l(x)}{z - \lambda_l} \int_{\partial\Omega} h(y) \, B\varphi_l(y) \, dS_y, \tag{16}$$

where $dS_y$ is the area element of the surface $\partial\Omega$ in $\mathbf{R}^m$ and, without loss of generality, we take real $\varphi_l$. However, the right hand side of (16) does not converge [NL]. Therefore, the sum in (16) has to be regularized. This was done by A. Nachman, J. Sylvester and G. Uhlmann in [NSU] for the case $g^{jk} = \delta^{jk}$, i.e. $a(x, D) = -\Delta + q$. Their construction was based on the fact that $-\Delta + q$ is a relatively compact perturbation of the Laplace operator $-\Delta$ and the boundary spectral data and the Dirichlet-to-Neumann maps for the Laplace operator in a given domain $\Omega$ can be always found. However, $-g^{1/2}\partial_j(g^{1/2}g^{jk}\partial_k) + q$ with unknown $g^{jk}$ is not a relatively compact perturbation of some known operator and thus one needs a different regularization. To this end, we start with the eigenfunction expansion

$$v_z^h = \sum_{l=1}^{\infty} \widehat{v}_l^h \varphi_l, \quad \widehat{v}_l^h = \frac{1}{z - \lambda_l} \langle h, B\varphi_l \rangle_{L^2(\partial\Omega)}, \tag{17}$$

which converges in $L^2(\Omega)$. Differentiating equation (17) with respect to $z$ we obtain a formula for the derivative $\partial_z v_z^h$,

$$\partial_z v_z^h = -\sum_{l=1}^{\infty} \frac{1}{(z - \lambda_l)^2} \langle h, B\varphi_l \rangle_{L^2(\partial\Omega)} \varphi_l. \tag{18}$$

Now the sum in the right hand side of (18) converges in $H^2(\Omega)$. Thus, the boundary spectral data determines the derivative of the Dirichlet-to-Neumann map,

$$\partial_z(\Lambda_z h) = -\sum_{l=1}^{\infty} \frac{1}{(z - \lambda_l)^2} \langle h, B\varphi_l \rangle_{L^2(\partial\Omega)} B\varphi_l|_{\partial\Omega}, \tag{19}$$

where the right hand side converges in $H^{1/2}(\partial\Omega)$. However, $\Lambda_z h$ can be represented as

$$\Lambda_z h = \lim_{\tau \to \infty} \left( \Lambda_{-\tau^2} h + \int_{-\tau^2}^{z} \partial_{z'}(\Lambda_{z'} h) \, dz' \right). \tag{20}$$

Thus, (19), (20) determine the map $\Lambda_z$ as soon as we can find the asymptotics of $\Lambda_{-\tau^2} h$, $\tau \to \infty$. It is shown in Appendix 1 that

$$(\Lambda_{-\tau^2} h)(x) = -\tau \rho(x) h(x) - \frac{1}{2}\rho(x) H(x) h(x) + O(\tau^{-1}), \quad \tau \to \infty, \tag{21}$$

where $\rho$ is given by formula (5) and $H(x)$ is the mean curvature of $\partial\Omega$ with respect to the metric $g_{jk}$. In turn, the mean curvature can be found from the symbol of $\partial_z \Lambda_z$ which is determined from the boundary spectral data by (19) (for a detailed construction see Appendix 1). Thus formulae (19)-(21) determine $\Lambda_z$ from the boundary spectral data.

<u>$ii. \to i.$</u> Formulae (18)-(19) show that $\Lambda_z$ is a meromorphic operator-valued function which has simple poles at the eigenvalues $z = \lambda_l$. This implies that $\Lambda_z$ determines the eigenvalues $\lambda_l$. Moreover, one can show that the residue of $\Lambda_z$ at $z = \lambda_l$ is a finite dimensional integral operator with the kernel

$$K_l(x, y) = \sum_{k \in L_l} B\varphi_k(x) B\varphi_k(y), \quad x, y \in \partial\Omega. \tag{22}$$

Here $L_l$ is the set of integers $k$ such that $\lambda_k = \lambda_l$. Clearly, when $\lambda_l$ has multiplicity one, the kernel $K_l(x, y)$ determines $B\varphi_l(x)$ up to a multiplication by $\pm 1$. In general, the linear independency of $B\varphi_k(x)$ makes it possible to find $B\varphi_k(x)$, $k \in L_l$ up to an orthogonal transformation. More precisely, we can find

$$\xi_j(x) = \sum_{k \in L_l} \alpha_{jk} B\varphi_k(x), \quad j \in L_l, \tag{23}$$

where $[\alpha_{jk}]$ is an orthogonal matrix (for details see [KKL]).

Clearly, $\xi_j$ are themselves the boundary values of some eigenfunctions $\psi_j$ which form another orthonormal basis of the eigenspace corresponding to the eigenvalue $z = \lambda_l$, i.e. $\xi_j = B\psi_j$. Thus, the Gel'fand data determines the boundary spectral data.

## 3 Equivalence of the hyperbolic and elliptic boundary data

<u>$ii \to iii.$</u> To construct $R^w$ from $\Lambda_z$, $z \in \mathbf{C}$, let

$$\widetilde{f}(x, \omega) = \int_0^\infty e^{-\omega t} f(x, t)\, dt \tag{24}$$

be the Laplace transform of $f \in C_0^\infty(\partial\Omega \times \mathbf{R}_+)$. It is well defined for $\omega \in \mathbf{C}$ and satisfies

$$||\widetilde{f}(\omega)||_{H^{3/2}(\partial\Omega)} \leq C_N (1 + |\omega|)^{-N}, \quad \operatorname{Re}\omega > 0, \tag{25}$$

for any $N > 0$. Let $v(\omega)$ be the solution of the Dirichlet problem (7) with $z = -\omega^2$ and $h = \widetilde{f}(\omega)$. Since

$$||v_z^h||_{H^2(\Omega)} \leq C(1 + |z|)\left(1 + \operatorname{dist}^{-1}(z, \sigma(A))\right) ||h||_{H^{3/2}(\partial\Omega)}, \tag{26}$$

estimate (25) implies that $v(\omega)$ decreases rapidly when $|\omega| \to \infty$ in the right half-plane, $\operatorname{Re}\omega \geq \mu > 0$. Moreover, $Bv(\omega)|_{\partial\Omega} = \Lambda_{-\omega^2}\widetilde{f}(\omega)$.

Using the inverse Laplace transform, we see that

$$u(x,t) = \frac{1}{2\pi} \int_{-\infty}^{\infty} e^{t(\mu+is)} v(x, \mu+is) ds$$

is well defined in $C^{\infty}(\mathbf{R}_+, H^2(\Omega))$ for $\mu > \sqrt{\mu_0}$, where

$$\mu_0 = \max\{0, -\lambda_1\}, \tag{27}$$

and solves initial boundary value problem (8). Therefore,

$$R^w f = Bu|_{\partial\Omega \times \mathbf{R}_+} = \frac{1}{2\pi} \int_{-\infty}^{\infty} e^{t(\mu+is)} \Lambda_{-(\mu+is)^2} \widetilde{f}(x, \mu+is) ds.$$

Thus $\Lambda_z$ determines the map $R^w$.

<u>$iii. \to iv.$</u> It follows from the definition of $\Pi^w$ that

$$\Pi^w f = \operatorname{Re} \int_{\partial\Omega} \int_0^{\infty} (R^w f)(x,t) \overline{\partial_t f(x,t)} \, dt dS_x, \tag{28}$$

which implies that the hyperbolic Dirichlet-to-Neumann map determines $\Pi^w$.

<u>$iv. \to i.$</u> The quadratic form $\Pi^w(f)$ determines the bilinear form

$$\Pi^w[f,h] = \frac{1}{2} \operatorname{Re} \int_{\partial\Omega} \int_0^{\infty} \left((R^w f)(x,t)\overline{\partial_t h(x,t)} + \partial_t f(x,t)\overline{(R^w h)(x,t)}\right) dt dS_x.$$

Because $R^w$ maps the real valued functions to the real valued functions, we can separate the real and imaginary parts of $\Pi^w$ and construct the complex bilinear form

$$\Pi^w_{\mathbf{C}}[f,h] = \frac{1}{2} \operatorname{Re} \int_{\partial\Omega} \int_0^{\infty} ((R^w f)(x,t)\partial_t h(x,t) + \partial_t f(x,t)(R^w h)(x,t)) \, dt dS_x. \tag{29}$$

Using the Plancherel formula we can represent the integral in the right-hand side of (29) in terms of the Fourier transforms $\widehat{f}, \widehat{h}$ of $f, h$, where e.g.

$$\widehat{f}(x,k) = \int_{\mathbf{R}} e^{ikt} f(x,t) \, dt.$$

Obviously, $\widehat{R^w f}(k) = \Lambda_{k^2}\widehat{f}(k)$. As $\Lambda_{k^2}$ has poles at points $k = \pm\sqrt{\lambda_l}$, the contour of integration in the Plancherel formula is shifted in the complex

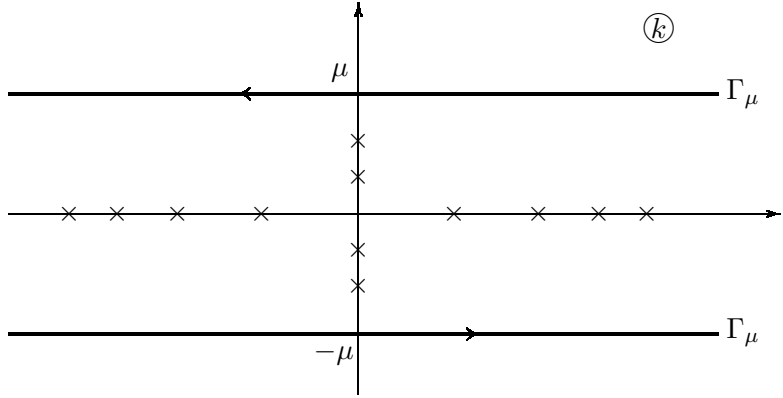
Figure 1: Contour $L_\mu$

domain to the boundary $\Gamma_\mu$ of the strip $\{k \in \mathbf{C} : |\operatorname{Im} k| \leq \mu\}$ (see Fig. 1) and we obtain

$$\Pi^w_{\mathbf{C}}[f,h] = \frac{1}{4\pi i} \int_{\Gamma_\mu} \int_{\partial\Omega} \Lambda_{k^2} \widehat{f}(-k)\,\widehat{h}(k)\,dS_x\,kdk \qquad (30)$$

(for details see [KKL]).

Taking into account the analyticity of $\widehat{f}(k)$ and $\widehat{h}(k)$ and the representation (22) of the residues of $\Lambda_z$ at $z = \lambda_l$, estimate (26) and Weyl's asymptotics for eigenvalues, the residue theorem yields that

$$\Pi^w_{\mathbf{C}}[f,h] = \frac{1}{2} \sum_{j=1}^{\infty} \left( \left\langle \widehat{f}(\sqrt{\lambda_j}), B\varphi_j \right\rangle \left\langle \widehat{h}(-\sqrt{\lambda_j}), B\varphi_j \right\rangle \right.$$
$$\left. + \left\langle \widehat{f}(-\sqrt{\lambda_j}), B\varphi_j \right\rangle \left\langle \widehat{h}(\sqrt{\lambda_j}), B\varphi_j \right\rangle \right). \qquad (31)$$

This formula provides a possibility to find the kernels $K_l(x,y)$ given by formula (22). Indeed, let $f = f_\tau$, $h$ be of the form

$$f_\tau(x,t) = F(x)\chi(t-\tau), \quad h(x,t) = H(x)\chi(t), \qquad (32)$$

where $F, H \in C^\infty(\partial\Omega)$ and $\chi \in C_0^\infty(\mathbf{R}_+)$ is real-valued. Then, formula (31) yields that

$$\Pi^w_{\mathbf{C}}[f_\tau, h] = \frac{1}{2} \sum_{l=1}^{\infty} \cos(\sqrt{\lambda_l}\,\tau)|\widehat{\chi}(\sqrt{\lambda_l})|^2 \langle K_l F, H \rangle_{L^2(\partial\Omega)}, \qquad (33)$$

where $K_l$ are the integral operators with kernels $K_l(x, y)$. Analyzing the asymptotic behavior of the sum (33) when $\tau \to \infty$, we can find $\lambda_l$ and $\langle K_l F, H \rangle$. As $F$ and $H$ are arbitrary, we find the kernels $K_l(x, y)$ for any $l$ (for details see [KKL]). Due to the implication $ii. \to i.$, this also determines the boundary spectral data.

## 4 Equivalence of the parabolic and non-stationary Schrödinger boundary data and the boundary spectral data

<u>$ii \to vi.$</u> Let $f \in C_0^\infty(\partial\Omega \times \mathbf{R}_+)$. Its Fourier transform with respect to time variable $t$, $\widehat{f}(x, k)$ is a rapidly decreasing function in the upper half plane $\mathbf{C}_+ = \{k \in \mathbf{C}: \operatorname{Im} k \geq 0\}$. Let

$$\eta(x, t) = \frac{1}{2\pi} \int_{i\mu + \mathbf{R}} e^{-ikt} \eta(x, k) \, dk, \quad \mu > 0,$$

where

$$a(x, D)\eta(x, k) = -k\eta(x, k) \quad \text{in } \Omega, \quad \eta(x, k)|_{\partial\Omega} = \widehat{f}(x, k).$$

Using estimate (26), we see that $\eta(x, t)$ is the solution of the initial boundary value problem (10), i.e., $\eta(x, t) = \psi^f(x, t)$. As $B\eta(k) = \Lambda_{-k}\widehat{f}(k)$, we obtain

$$R^s f = \frac{1}{2\pi} \int_{i\mu + \mathbf{R}} e^{-ikt} \Lambda_{-k} \widehat{f}(k) \, dk. \tag{34}$$

<u>$vi \to vii.$</u> By direct calculations we obtain

$$\Pi^s(f) = \operatorname{Re} \int_{\partial\Omega} \int_0^\infty R^s f(x, t) \overline{\partial_t f(x, t)} \, dt dS_x. \tag{35}$$

Thus, $R^s$ determines $\Pi^s$.

For the future we note that $R^s$ determines also

$$\|\psi^f(\infty)\|_{L^2} = \lim_{t \to \infty} \|\psi^f(t)\|_{L^2}.$$

Indeed

$$\|\psi^f(\infty)\|_{L^2} = i \int_{\partial\Omega} \int_0^\infty \left( R^s f(x, t) \overline{f(x, t)} - f(x, t) \overline{R^s f(x, t)} \right) dt dS =$$

$$2 \operatorname{Im} \int_{\partial\Omega} \int_0^\infty f(x,t) \overline{R^s f(x,t)}\, dt dS. \tag{36}$$

This implies that $R^s$ determines also the energy flux $\Pi^s_{E_0}(f)$ which corresponds to the potential $q + E_0$, $E_0 \in \mathbf{R}$.

<u>$vii \to i.$</u> Quadratic form $\Pi^s(f)$ determines the form $\Pi^s_r[f,h]$,

$$\Pi^s_r[f,h] = \frac{1}{2}\operatorname{Re} \int_{\partial\Omega} \int_0^\infty \left(R^s f(x,t)\overline{\partial_t h(x,t)} + \overline{\partial_t f(x,t)}\, R^s h(x,t)\right) dt dS_x =$$

$$\frac{1}{2}\operatorname{Re} \int_{\partial\Omega} \int_0^\infty \left(R^s f(x,t)\overline{\partial_t h(x,t)} + \partial_t f(x,t)\, \overline{R^s h(x,t)}\right) dt dS_x.$$

Taking $\widetilde{f} = if$ in this formula we obtain also that

$$\Pi^s_r[\widetilde{f},h] = -\frac{1}{2}\operatorname{Im} \int_{\partial\Omega} \int_0^\infty \left(R^s f(x,t)\overline{\partial_t h(x,t)} + \partial_t f(x,t)\, \overline{R^s h(x,t)}\right) dt dS_x$$

Thus, the quadratic form $\Pi^s(f)$ determines the bilinear form $\Pi^s[f,h]$,

$$\Pi^s[f,h] = \frac{1}{2} \int_{\partial\Omega} \int_0^\infty \left(R^s f(x,t)\overline{\partial_t h(x,t)} + \partial_t f(x,t)\, \overline{R^s h(x,t)}\right) dt dS_x \tag{37}$$

$$= \frac{1}{2} \int_{\partial\Omega} \int_0^\infty \left(e^{-\mu t} R^s f(x,t) \overline{e^{\mu t}\partial_t h(x,t)} + e^{\mu t}\partial_t f(x,t)\, \overline{e^{-\mu t} R^s h(x,t)}\right) dt dS_x,$$

where $\mu > 0$ is arbitrary. Then, by the Parseval identity,

$$\Pi^s[f,h] = \frac{1}{2} \int_{\partial\Omega} \int_{-\infty}^\infty \left((e^{-\mu t}R^s f)^\wedge(x,k)\overline{(e^{\mu t}\partial_t h)^\wedge(x,k)} + \right.$$
$$\left. + (e^{\mu t}\partial_t f)^\wedge(x,k)\overline{(e^{-\mu t}R^s h)^\wedge(x,k)}\right) dk dS_x.$$

Using formula (34), we see that

$$\Pi^s[f,h] = \frac{1}{2} \int_{-\infty}^\infty \left\{-i(k-i\mu)\langle \widehat{f}(k-i\mu), (\Lambda_{-k-i\mu}\widehat{h})(k+i\mu)\rangle \right.$$

$$\left. + i(k+i\mu)\,\langle (\Lambda_{-k-i\mu}\widehat{f})(k+i\mu), \widehat{h}(k-i\mu)\rangle\right\} dk =$$

$$= \frac{1}{4\pi i} \int_{\partial\Omega} \int_{\Gamma_\mu} z\, \Lambda_{-z}\widehat{f}(x,z)\overline{\widehat{h}(x,\overline{z})}\, dz dS_x, \tag{38}$$

where $\Gamma_\mu$ is given in section 3. In deriving the last equation in (38) we use the identity $(\Lambda_{-z})^* = \Lambda_{-\overline{z}}$.

Since $\Lambda_{-z}$ has poles at points $z = -\lambda_l$, we can apply the residue theorem as in section 3 and obtain

$$\Pi^s[f, h] = \frac{1}{2} \sum_{l=1}^{\infty} \lambda_l \langle K_l \widehat{f}(-\lambda_l), \widehat{h}(-\lambda_l) \rangle.$$

Taking $f$ and $h$ of form (32) with $F, G \in C^\infty(\partial\Omega)$, we get

$$\Pi^s[f, h] = \frac{1}{2} \sum_{l=1}^{\infty} \lambda_l \exp(-i\lambda_l \tau) \langle K_l F, H \rangle |\widehat{\chi}(\lambda_l)|^2. \tag{39}$$

The sum (39) considered as a function of $\tau \in \mathbf{R}$ determines all the eigenvalues $\lambda_l$, $\lambda_l \neq 0$ and the corresponding bilinear forms $K_l$ (see section 3 for further details). Thus $\Pi^s$ determines the boundary spectral data except for $B\varphi_l$ with $\lambda_l = 0$.

We note that the terms $B\varphi_l$ for which $\lambda_l = 0$ do not appear in series (39) and thus the boundary spectral data corresponding to the eigenvalue 0 can not be determined from the asymptotics.

<u>$vi. \to i.$</u> As noted after formula (36), having $R^s$ it is possible to construct the form $\Pi^s_{E_0}$. Thus, by choosing sufficiently large $E_0$ and using formulae (35), (36) we can determine all the eigenvalues $\lambda_l$ and corresponding functions $B\varphi_l$. Similarly if, in addition to the energy flux $\Pi^s$, we know the quadratic form $f \mapsto \|\psi^f(\infty)\|^2$, we can also determine all $\lambda_l$, $B\varphi_l$.

<u>$ii. \to v.$</u> Using the Laplace transform instead of the Fourier transform and following the same steps as for the non-stationary Schrödinger operator, we obtain

$$R^h f(x, t) = \frac{1}{2\pi i} \int_{\mu + i\mathbf{R}} e^{\omega t} \Lambda_{-\omega} \widetilde{f}(\omega) \, d\omega. \tag{40}$$

(compare with formula (34)).

<u>$v. \to i.$</u> Similar considerations as for the non-stationary Schrödinger operator show that $\widetilde{R^h f}(x, \omega) = \Lambda_{-\omega} \widetilde{f}(x, \omega)$ for $\operatorname{Re}\omega > \mu_0$ where $\mu_0$ is of form (27). Using, as earlier, the boundary sources of form (32) and applying the Parseval formula and the residue theorem, we obtain

$$\int_{\partial\Omega} \int_0^\infty f(x, t) R^h h(x, t) \, dS_x dt = \sum_{l=1}^\infty e^{-\lambda_l \tau} |\widetilde{\chi}(\lambda_l)|^2 \langle K_l F, H \rangle \tag{41}$$

(for details, see [KKL]). As above, the asymptotics of the sum in the right hand side of (41) for $\tau \to \infty$ determines the boundary spectral data. This proves Theorem 1.

## 5 Generalizations

### 5.1 Continuation of data.

From the point of view of uniqueness, we are often interested in the case when the data is given only on a finite time or spectral intervals. In the other words, we are given only the restrictions $R^w f|_{\partial\Omega \times [0,T]}$, $R^s f|_{\partial\Omega \times [0,T]}$ or $R^h f|_{\partial\Omega \times [0,T]}$ on the time interval $[0,T]$ for the time domain problems or $\Lambda_z$ for $z \in I$ where $I \subset \mathbf{C}$ for the spectral problem. In particular, as inverse problems are usually ill-posed, it is important to find out whether these restricted data determine the complete data, i.e. $R^w$, $R^s$, $R^h$ or $\Lambda_z$ for $z \in \mathbf{C}$.

In the spectral problem we see that $\Lambda_z$ is a meromorphic operator-valued function of $z$. Thus, by analytic continuation $\Lambda_z$, $z \in I$ determines $\Lambda_z$ for all $z \in \mathbf{C}$ as soon as $I$ has an accumulation point. For the time-domain problems, we can use a complexification of time for the similar purposes. Indeed, when $f \in C_0^\infty(\partial\Omega \times [0, T_0])$, we can use the Paley-Wiener theorem and estimate (26) to deform the countours of integration in integrals (34) and (40) with $t > T_0$ to the contour $\widetilde{\Gamma_\mu}$, $\mu > \mu_0$. This contour consists of a straight ray $(-\infty - i\mu, -i\mu)$, semicircle of the radius $\mu$ and straight ray $(i\mu, -\infty + i\mu)$. Then for any function $f \in C_0^\infty(\partial\Omega \times [0, T_0])$, the functions

$$t \mapsto R^h f(x,t), \quad t \in \{\operatorname{Re} t > T_0\},$$
$$t \mapsto R^s f(x,t), \quad t \in \{\operatorname{Re} t > T_0, \ \operatorname{Im} t > 0\}$$

are holomorphic functions. Moreover, $R^s f(x,t)$ is continuous with respect to $t \in \mathbf{C}$ up to the real semiaxis $t > T$. Thus, for any $f \in C_0^\infty(\partial\Omega \times [0, T_0])$ $R^h f|_{\partial\Omega \times [0,T]}$ and $R^s f|_{\partial\Omega \times [0,T]}$, $T > T_0$ determine, correspondingly, $R^h f$ and $R^s f$ on $\partial\Omega \times \mathbf{R}_+$. However, $R^s$ and $R^h$ are linear operator which commute with the time delay operator $Y_\tau$, $Y_\tau f(t) = f(t - \tau)$. Thus, the restrictions $R^h f|_{\partial\Omega \times [0,T]}$ and $R^s f|_{\partial\Omega \times [0,T]}$ with $f \in C_0^\infty(\partial\Omega \times [0, T_0])$ determine $R^h f$ and $R^s f$ for any $f \in C_0^\infty(\partial\Omega \times \mathbf{R}_+)$.

In the case of the wave equation the situation is more complicated. However, it is known that given $R^w f|_{[0,T]}$ for any $f \in C_0^\infty(\partial\Omega \times [0,T])$ we can construct $R^w f|_{\mathbf{R}_+}$ for all $f \in C_0^\infty(\partial\Omega \times \mathbf{R}_+)$ as soon as $T > 2 \max\{d(x, \partial\Omega) : x \in \Omega\}$ (see [KL2], [KKL]). We note that in this case it is also possible to solve the inverse problem for the wave equation directly, without continuation of data. For this fact we refer to e.g. [BKa], [B], [KKL], and, for a more general case of a wave equation which corresponds to a non-selfadjoint symbol $a(x, D)$, to [KL1].

## 5.2 Anisotropic boundary form and energy flux on a finite time interval.

Let us consider the inverse problem for the wave equation with finite time observations. It is known that with appropriate boundary data measured on a time-interval $[0, T]$ it is possible to reconstruct the operator in the domain $\{x : d(x, \partial\Omega) < T/2\}$. In the reconstruction procedure it is actually not necessary to know the Dirichlet-to-Neumann map $R^w \widetilde{f}|_{[0,T]}$ but only the antisymmetric form corresponding to it (see [KKL]). This antisymmetric form, $B^T[f, h]$ is given by

$$B^T[f,h] = \int_{\partial\Omega} \int_0^T (f(x,t)R^w h(x,t) - h(x,t)R^w f(x,t)) \, dt dS_x$$

with real-valued $f, h \in C^{0,\infty}(\partial\Omega \times \mathbf{R}_+)$, where $C^{0,\infty}(\partial\Omega \times \mathbf{R}_+)$ consists of $C^\infty$-functions with support in $\mathbf{R}_+$. As this antisymmetric form is related to the symplectic structure on $L^2(\partial\Omega \times [0,T])$ it is often a more natural object in inverse problems than the Dirichlet-to-Neumann map.

Next we consider the relations between $B^T$ and the finite time energy flux,

$$\Pi^T[f,h] = \frac{1}{2} \operatorname{Re} \int_{\partial\Omega} \int_0^T (R^w f(x,t)\partial_t h(x,t) + \partial_t f(x,t) R^w h(x,t)) \, dt dS_x \quad (42)$$

also considered on real-valued $f$ and $h$. Hence, $\Pi^T[f, f] = E^w(u^f, T)$ is the amount of energy that has passed into $\Omega$ through its boundary until $t = T$. In particular, we consider the question of finding the antisymmetric boundary form $B^T[f, h]$ from the energy flux $\Pi^T[f, h]$ of form (42) given for $f, h \in C^{0,\infty}(\partial\Omega \times \mathbf{R}_+)$.

For $f, h \in C^{0,\infty}(\partial\Omega \times \mathbf{R}_+)$, $\partial_t R^w f = R^w \partial_t f$ so that

$$B^T[\partial_t f, h] = \int_{\partial\Omega} \int_0^T \partial_t f(x,t) R^w h(x,t) \, dt dS_x + \quad (43)$$

$$+ \int_{\partial\Omega} \int_0^T \partial_t h(x,t) R^w f(x,t) \, dt dS_x - \int_{\partial\Omega} h(x,T) R^w f(x,T) \, dS_x =$$

$$= 2\Pi^T[f,h] - \int_{\partial\Omega} h(x,T) R^w f(x,T) \, dS_x.$$

As any $f \in C^{0,\infty}(\partial\Omega \times \mathbf{R}_+)$ is of form $f = \partial_t f_0$, $f_0 \in C^{0,\infty}(\partial\Omega \times \mathbf{R}_+)$, representation (43) makes it possible to evaluate $B^T[f,h]$ in terms of $\Pi^T[f_0, h]$ as soon $h(x,T) = 0$. Similar considerations show that $\Pi^T$ determines also $B^T[f,h]$ when $f(x,T) = 0$. We need also the following result proven in Appendix 2.

**Proposition 1** Let $F, H \in C^\infty(\partial\Omega)$ be real-valued and $\alpha > 2$. Then

$$B^T[t^\alpha F, t^\alpha H] = \frac{1}{\alpha+1} \left( \Pi^T[t^\alpha F, t^{\alpha+1} H] - \Pi^T[t^{\alpha+1} F, t^\alpha H] \right). \tag{44}$$

Because any $f \in C^{0,\infty}(\partial\Omega \times \mathbf{R}_+)$ can be represented as

$$f(x,t) = f_0(x,t) + t^\alpha F(x), \quad f_0(x,T) = 0,\ \alpha > 2,$$

and similar representation takes place for $h$, we see:

**Proposition 2** Knowing $\Pi^T[f,h]$ for all $f, h \in C^{0,\infty}(\partial\Omega \times \mathbf{R}_+)$ we can determine $B^T[f,h]$.

### 5.3 Riemannian manifolds.

We mention that all previous results can be directly generalized to manifolds. Indeed, let $(M, g)$ be a Riemannian manifold with $-\Delta_g$ being its Laplace-Beltrami operator. The corresponding Schrödinger operators on the manifold $M$ are operators of the form $A = -\Delta_g + q$. Then the previous arguments remain valid for Schrödinger operators on $M$. In particular, Theorem 1 claims that all data *i.-vi.* (and also *vii.* when zero is not an eigenvalue) are equivalent.

### 5.4 Gauge transformations

When we study equivalence of boundary data for general 2nd-order operators instead of the Schrödinger operators (1)-(2), we have to take into account possible coordinate and gauge transformations. A general theory of gauge transformations in inverse problems for an arbitrary self-adjoint 2nd-order operator

$$Au = \sum_{j,k=1}^m a^{jk}(x)\partial_j\partial_k u + \sum_{j=1}^m b^j(x)\partial_j u + c(x)u, \tag{45}$$

$$\mathcal{D}(\widetilde{A}) = H^2(\Omega) \cap H^1_0(\Omega)$$

is developed in [KKL], see also [K1], [KK]. In this subsection we will consider a subclass of operators (45) and show how the results obtained for Schrödinger operators can be generalized for them. Naturally, when studying a subclass of operators (45), we have to consider the subgroup of the complete group of the coordinate and gauge transformations which preserves this subclass. This subgroup is called the admissible group.

As an example, we consider the class of the anisotropic conductivity operators,

$$a(x, D)v(x) = -\sum_{j,k=1}^{m} \partial_j \left(a^{jk}(x)\partial_k v(x)\right) + q(x)v(x) \qquad (46)$$

with the Dirichlet boundary condition (2), rather than a Schrödinger operator (1). We note that all operators of form (45) that are self-adjoint with respect to the Lebesque measure on $\Omega$ can be written in form (46). Then

$$\{\lambda_l, \ B\varphi_l|_{\partial\Omega}: \ l=1,2,\dots\}, \quad B\varphi_l = \sum_{j,k=1}^{m} a^{jk} n_j \partial_k \varphi_l, \qquad (47)$$

is the corresponding boundary spectral data. Similarly, we can define elliptic, hyperbolic, etc. boundary data of form $ii.$-$vii.$ with the boundary operator $B$ of form (47). Although it seems that the question of the equivalence between data i.-vii. for operators (46), (47) is completely analogous to that for Schrödinger operators, (1), (2) there is a significant difference. Namely, it is, in principle, impossible to determine the corresponding Dirichlet-to-Neumann map,

$$f \to Bu_z^f|_{\partial M}$$

with $B$ of form (47) and also the corresponding non-stationary Dirichlet-to-Neumann maps $R^w$, $R^h$, and $R^s$ from the boundary spectral data (47). As we will see below, the reason for this is that the admissible group of gauge and coordinate transformations which does not change the type of an operator under consideration can change the boundary data. More precisely, let $X: \Omega \to \Omega$, $X|_{\partial\Omega} = id|_{\partial\Omega}$ be a diffeomorphism (coordinate transformation) in the domain $\Omega$. It gives rise to the transformation, $S_X: L^2(\Omega) \to L^2(\Omega)$,

$$S_X u(x) = u(X(x)).$$

Similarly, any $\kappa \in C^\infty(\Omega)$, $\kappa > 0$ in $\Omega$, gives rise to the gauge transformation, $S_\kappa: L^2(\Omega) \to L^2(\Omega)$,

$$S_\kappa u(x) = \kappa(x) u(x).$$

Then, $S_\kappa \circ S_X : L^2(\Omega) \to L^2(\Omega)$ induces a transformation for operators which maps $a(x, D)$ of form (45) to the operator

$$\widetilde{a}(x, D)u = S_\kappa \circ S_X \, a(x, D) \left( S_X^{-1} \circ S_\kappa^{-1} u \right). \tag{48}$$

Next, assume that $a(x, D)$ is of form (46). Taking

$$\kappa(X(x)) = |\det(dX)|^{-1}(x), \tag{49}$$

we see that $\widetilde{a}(x, D)$ is of the same form (46) as $a(x, D)$ with, however, different coefficients $\widetilde{a_{ij}} \neq a_{ij}$, $\widetilde{q} \neq q$. Furthermore, it can be shown that the only transformations of the form $S_\kappa \circ S_X$ which preserve the form (46) of an operator in a given $\Omega$ are given by (48)-(49). Moreover, if $\det(dX) = 1$ on $\partial\Omega$, the boundary spectral data of $a(x, D)$ and $\widetilde{a}(x, D)$ coinside. On the other hand, in this case

$$\widetilde{\Lambda_z} f = \Lambda_z f + (B\kappa)|_{\partial\Omega} f,$$

Analysing the proof of Lemma 1 in Appendix 1, we see that asymptotics (21) remains valid for the operator (46) with, however, other factors $\alpha(x)$ and $K(x)$ instead of $\rho(x)$ and $H(x)$. These factors depend only on $a_{ij}(x)$. We can then follow the arguments of section 2, $i. \to ii.$ to show that the boundary spectral data determine the Dirichlet-to-Robin map

$$f \to \widetilde{B} u_z^f |_{\partial\Omega} = B u_z^f + \sigma u_z^f |_{\partial\Omega}, \, z \in \mathbf{C}$$

with an unknown function $\sigma \in C^\infty(\partial M)$. Motivated by this we say that operators $\Lambda_z^{(1)}$ and $\Lambda_z^{(2)}$ are gauge-equivalent if there is a function $\sigma$ such that $\Lambda_z^{(1)} f = \Lambda_z^{(2)} f + \sigma f$. Similarly, two bilinear forms $\Pi_z^{(1)}$ and $\Pi_z^{(2)}$ are gauge-equivalent if there is a function $\sigma$ such that

$$\Pi_z^{(1)}(f, f) = \Pi_z^{(2)}(f, f) + \int_{\partial\Omega} \sigma(x)|f(x)|^2 \, dS_x.$$

Then, using similar considerations as in sections 2-4 we obtain:

**Proposition 3** *Inverse problems i.-vi. for anisotropic conductivity operators (46) are gauge-equivalent, in the sense that the data i. determines the gauge-equivalence class of any data i.-vi. and vise versa.*

We would note that this phenomenon does not occur in the Schrödinger case because its admmissible group of transformations consists of only changes of coordinates, that is, $S_X$.

In connection to the inverse problem, we would mention that it is possible to find $a(x, D)$ (upto a transformation from the group (48)) from either the boundary spectral data (47) or the Dirichlet-to-Robin map $\Lambda_z + \sigma(x)I$, $z \in \mathbf{C}$ with unknown $\sigma$ (see e.g. [K1], [KK], [KKL]).

# 6 Appendix 1

In this appendix we will prove asymptotic expansion (21). To this end we introduce boundary normal coordinates $(y, n)$, $y = (y^1, \cdots, y^{m-1})$, $n \geq 0$ in a vicinity of $\partial \Omega$. Let $x \in \Omega$. Denote by $n = n(x)$ the distance (in the metric $g$) from $x$ to the boundary $\partial \Omega$, $n(x) = d_g(x, \partial \Omega)$. When $x$ is sufficiently close to $\partial \Omega$ there is a unique point $y = y(x)$ with $n(x) = d_g(x, y(x))$. Introducing some (local) coordinates $(y^1, \cdots, y^{m-1})$ on $\partial \Omega$ we use $(y^1(x), \cdots, y^{m-1}(x), n(x))$ as (boundary normal) coordinates of $x$. The length element $dl$ in these coordinates has a special form,

$$dl^2 = dn^2 + g_{\alpha\beta}(y, n)\, dy^\alpha dy^\beta, \quad \alpha, \beta = 1, \cdots, m-1.$$

Then

$$H(y, n_0) = \frac{1}{2} \partial_n \left( \log g(y, n) \right) |_{n=n_0} \tag{50}$$

is the mean curvature of the surface $\partial \Omega_{n_0} = \{x \in \Omega : n(x) = n_0\}$. In particular, $\partial \Omega_0 = \partial \Omega$.

Let $a(x, D)$ be of the form

$$a(x, D)v = -\Delta_g v + q(x)v.$$

We denote by $D_y = -i\partial_y$, $D_n = -i\partial_n$, $D_s = -i\partial_s$ the derivatives in boundary normal coordinates. The crucial step in proving (21) is the following lemma.

**Lemma 1** *The Dirichlet-to-Neumann map $\Lambda_z$ of form (4), (6), (7) for $a(x, D)$ has an asymptotic expansion*

$$\Lambda_z \sim \sum_{k=0}^{\infty} \mathcal{P}_k(y, D_y)(-z)^{-(k-1)/2}, \tag{51}$$

*where $\mathcal{P}_k(y, D_y)$ are differential operators of order $k$. Expansion (51) means that*

$$||\Lambda_z - \sum_{k=0}^{l} \mathcal{P}_k(y, D_y)(-z)^{-(k-1)/2}||_{H^{l+1}(\partial\Omega) \to L^2(\partial\Omega)} \leq C_l |z|^{-l/2}$$

*in the domain $|\arg(-z)| \leq \delta\pi$, $0 < \delta < 1$. In particular,*

$$\mathcal{P}_0(y, D_y) = -\rho, \quad \mathcal{P}_1(y, D_y) = -\frac{H(y)\rho}{2}, \tag{52}$$

$$\mathcal{P}_2(y, D_y) = -\frac{1}{2}\rho g^{\alpha\beta}(y, 0) D_\alpha D_\beta + \widetilde{\mathcal{P}}_2(y, D_y),$$

$$\mathcal{P}_3(y, D_y) = -\frac{1}{4}\rho \partial_n g^{\alpha\beta}(y, 0) D_\alpha D_\beta + \widetilde{\mathcal{P}}_3(y, D_y),$$

where $H(y) = H(y, 0)$ and $\widetilde{\mathcal{P}}_2$, $\widetilde{\mathcal{P}}_3$ are differential operators of order 1 with coefficients depending on the derivatives of the metric tensor and the potential on $\partial\Omega$.

**Proof.** By adding, if necessary, a positive constant to $q(x)$ we can assume that $\sigma(A) \subset \{z : \operatorname{Re} z > 0\}$ where $A$ the operator with symbol $a(x, D)$ and the Dirichlet boundary condition.

In [LU] it is shown that the Dirichlet-to-Neumann map for an elliptic differential operator of the second order is a pseudodifferential operator. Its symbol depends on the derivatives of the metric tensor at the boundary. To apply this approach, we introduce an auxiliary elliptic operator $\widetilde{a}(x, s, D_x, D_s) = \widetilde{a}(x, D_x, D_s)$,

$$\widetilde{a}(x, D_x, D_s) = a(x, D_x) - \partial_s^2 \quad \text{in} \quad (x, s) \in \widetilde{\Omega} = \Omega \times \mathbf{R}.$$

Rewriting $\widetilde{a}(x, D_x, D_s)$ in coordinates $(y, s, n)$ we have

$$\begin{aligned}\widetilde{a}(x, D_x, D_s) &= \widetilde{a}(y, n, D_y, D_n, D_s) \\ &= -\partial_n^2 - \partial_s^2 - g^{\alpha\beta}(y, n)\partial_\alpha \partial_\beta - H(y, n)\partial_n + f^\alpha(y, n)\partial_\alpha + h(y, n).\end{aligned}$$

By [LU], $\widetilde{a}(x, s, D_x, D_s)$ has a factorization

$$\widetilde{a} = (D_n - iH(y, n) - iL(y, n, D_y, D_s))(D_n + iL(y, n, D_y, D_s)), \qquad (53)$$

The operator $L(y, n, D_y, D_s)$ is a pseudodifferential operator of the first order with respect to $(y, s)$ which depends smoothly on the parameter $n$. The symbol of $L(y, n, D_y, D_s)$ has the form

$$L(y, n, \xi, \tau) \sim \sum_{m \leq 1} L_m(y, n, \xi, \tau), \quad y \in \mathbf{R}^{m-1}, \ n \in \mathbf{R}_+,$$

where $\xi \in \mathbf{R}^{m-1}$, $\tau \in \mathbf{R}$ are dual to $y$ and $s$, correspondingly. $L_m(y, s, \xi, \tau)$, $m = 1, 0, -1, \ldots$, are positive-homogeneous symbols of order $m$ with respect to $(\xi, \tau)$ and we take the principal symbol $L_1(y, n, \xi, \tau)$ to be negative.

Using the calculus of pseudodifferential operators and comparing the terms of the same homogeneity with respect to $(y, \tau)$ on the both sides of equality (53), we obtain

$$\sum_{\substack{j,k \leq 1 \\ j+k-|\gamma|=m}} \frac{1}{\gamma!} \partial_\xi^\gamma L_j D_y^\gamma L_k + (\partial_n + H)L_m = \widetilde{a}_m, \qquad (54)$$

where $\widetilde{a}_m = \widetilde{a}_m(y, n, D_y, D_s)$ is the term of homogeneity $m$ with respect to $(D_y, D_s)$ in $\widetilde{a}$. In particular, $\widetilde{a}_m = 0$ for $m < 0$. Using equation (54) we can find $L_m(y, n, \xi, \tau)$ in terms of coefficients of operator $\widetilde{a}$. In particular, equation (54) for $m = 2$ implies that

$$L_1(y, n, \xi, \tau) = -\widetilde{\eta} = -\sqrt{\eta^2 + \tau^2}, \quad \eta = \sqrt{g^{\alpha\beta}\xi_\alpha\xi_\beta}. \tag{55}$$

Further equations (54) give us recurrently that

$$L_m = \frac{1}{2\widetilde{\eta}} \left( \sum_{\substack{m+1 \leq j, k \leq 1 \\ j+k-|\gamma|=m+1}} \frac{1}{\gamma!} \partial_\xi^\gamma L_j D_y^\gamma L_k + (\partial_n + H) L_{m+1} - \widetilde{a}_{m+1} \right) \tag{56}$$

Formula (56) implies that

$$L_m(y, n, \xi, \tau) = \sum_{j_0(m) \leq j \leq m} p^{(m)}_{m-j}(y, n, \xi) \widetilde{\eta}^j, \quad m \leq 0, \tag{57}$$

where $p^{(m)}_k(y, n, \xi)$ are polynomials of order $k$ with respect to $\xi$. Using the binomial Taylor expansion of $\widetilde{\eta}$ when $\tau \to \infty$,

$$\widetilde{\eta}^k = \tau^k \sum_{j=0}^m \binom{k/2}{j} \tau^{-2j} \eta^{2j} + \mathcal{O}(\tau^{k-2m-2} \eta^{2m+2}), \tag{58}$$

where $\binom{a}{k} = a(a-1)\ldots(a-k+1)/k!$, we see that

$$L_m(y, n, \xi, \tau) = \sum_{j \leq m} \widehat{p}^{(m)}_{m-j}(y, n, \xi) \tau^j. \tag{59}$$

Here $\widehat{p}^{(m)}_k(y, s, \xi)$ are again polynomials of order $k$ with respect to $\xi$. These polynomials can be found by direct computations. In particular, from form (55) for $L_1$ we see that

$$\widehat{p}^{(1)}_1(y, n, \xi) = -1, \quad \widehat{p}^{(1)}_{-1}(y, n, \xi) = -\frac{1}{2} g^{\alpha\beta}(y, n) \xi_\alpha \xi_\beta. \tag{60}$$

Furthermore, analysing formula (56) with $m = 0$ we obtain that

$$\widehat{p}^{(0)}_0(y, n, \xi) = -\frac{H(y, n)}{2}, \quad \widehat{p}^{(0)}_{-2}(y, n, \xi) = -\frac{1}{4} \partial_n g^{\alpha\beta}(y, n) \xi_\alpha \xi_\beta. \tag{61}$$

Since for $L_1$ with the negative principal symbol,
$$L(y, 0, D_y, D_s) = \widetilde{\Lambda_0},$$
where $\widetilde{\Lambda_z}$ is the Dirichlet-to-Neumann map for $\widetilde{a}$, i.e.
$$\widetilde{\Lambda_z} \widetilde{f} = \partial_n \widetilde{v_z^f}|_{\partial\Omega}, \quad \text{where } \widetilde{a}(x, D)\widetilde{v_z^f} = z\widetilde{v_z^f}, \ \widetilde{v_z^f}|_{\partial\Omega} = \widetilde{f},$$
formulae (57), (59) (with $n = 0$) determine the symbol of $\widetilde{\Lambda_0}$.

Now we are ready to derive the asymptotics of the Dirichlet-to-Neumann map $\Lambda_{-\tau^2}$ of the operator $A$ when $\tau \to \infty$, $\tau \in \mathbf{R}$. Consider the boundary value problem
$$\begin{aligned} \widetilde{a}(x, D_x, D_s)\widetilde{v} &= 0, \quad \text{in } \widetilde{\Omega}, \\ \widetilde{v}|_{\partial\widetilde{\Omega}} &= f(y)e^{i\tau s}. \end{aligned}$$
Its periodic solution $\widetilde{v}(x, s) = v(x)e^{i\tau s}$ has the function $v(x)$ being the solution of the Dirichlet problem
$$\begin{aligned} a(x, D)v + \tau^2 v &= 0, \quad \text{in } \Omega, \\ v|_{\partial\Omega} &= f. \end{aligned}$$
Since $Bv = \rho \frac{\partial v}{\partial n}\big|_{n=0}$, we have
$$\Lambda_{-\tau^2} f = \rho e^{-i\tau s} \widetilde{\Lambda_0}(e^{i\tau s} f). \tag{62}$$
Therefore, $\Lambda_{-\tau^2}$ is a pseudodifferential operator with symbol
$$\lambda_{-\tau^2}(y, \xi) = \rho(y) L(y, 0, \xi, \tau). \tag{63}$$
Summarizing the previous considerations, we obtain an asymptotic expansion for the operator $\Lambda_{-\tau^2}$,
$$\Lambda_{-\tau^2} = \sum_{k=0}^{l} \mathcal{P}_k(y, D_y) \tau^{-k+1} + E_l(y, D_y; \tau), \tag{64}$$
where $\mathcal{P}_k(y, D_y)$ is a differential operator of order $k$ and $E_l(y, D_y; \tau)$ is a bounded operator from $H^{l+1}(\partial\Omega)$ to $L^2(\partial\Omega)$ with the norm $||E_l|| \leq C_l \tau^{-l}$. The exact form of the operators $\mathcal{P}_k(y, D_y)$ can be obtained by direct calculations using formulae (55)-(59), (62), (62). In particular, the first four polynomials have form (52). □

Lemma 1 yields immediately the following result.

**Corollary 1** *The derivative of Dirichlet-to-Neumann mapping, $\partial_z \Lambda_z$ has the expansion*

$$\partial_z \Lambda_z \sim \frac{1}{2}\rho(-z)^{-1/2} + \sum_{k=2}^{\infty} \frac{k-1}{2} \mathcal{P}_k(y, D_y)(-z)^{-(k+1)/2}, \qquad (65)$$

*where $P_k(y, D_y)$ are differential operators given in Lemma 1. In particular, the map $\partial_z \Lambda_z$ determines the weight $\rho$ and mean curvature $H(y)$.*

**Proof.** From the principal symbol of $\partial_z \Lambda_z$ we obtain $\rho$. It is then clear from formula (52) for $\mathcal{P}_3(y, D)$ that the map $\partial_z \Lambda_z$ determines the mean curvature $H(y)$. $\square$

# 7 Appendix 2

In this appendix we give the proof of Proposition 1. We can represent $B^T[t^\alpha F, t^\alpha H]$ in the form

$$B^T[t^\alpha F, t^\alpha H] = -\int_0^T \int_{\partial\Omega} t^\alpha \left(\partial_\nu u^{t^\alpha F} H - F \partial_\nu u^{t^\alpha H}\right) dt dS_x =$$

$$\frac{1}{\alpha+1} \int_0^T \int_{\partial\Omega} \left(\partial_\nu u^{t^{\alpha+1}F} \partial_t(t^\alpha H) + \partial_t(t^{\alpha+1}F) \partial_\nu u^{t^\alpha H}\right) dt dS_x -$$

$$- \frac{T^\alpha}{\alpha+1} \int_{\partial\Omega} \partial_\nu u^{t^{\alpha+1}F}(x, T) H(x) dS_x. \qquad (66)$$

Similarly,

$$B^T[t^\alpha F, t^\alpha H] = -\frac{1}{\alpha+1} \int_0^T \int_{\partial\Omega} \left(\partial_\nu u^{t^{\alpha+1}H} \partial_t(t^\alpha F) + \partial_t(t^{\alpha+1}H) \partial_\nu u^{t^\alpha F}\right) dt dS_x +$$

$$+ \frac{T^\alpha}{\alpha+1} \int_{\partial\Omega} \partial_\nu u^{t^{\alpha+1}H}(x, T) F(x) dS_x. \qquad (67)$$

Subtracting (67) from (66) and using definition (42) of $\Pi^T$, we obtain that

$$\int_{\partial\Omega} \left(\partial_\nu u^{t^{\alpha+1}F}(x, T) H(x) + \partial_\nu u^{t^{\alpha+1}H}(x, T) F(x)\right) dS_x =$$

$$2T^{-\alpha}\left(\Pi^T[t^{\alpha+1}F, t^\alpha H] + \Pi^T[t^\alpha F, t^{\alpha+1}H]\right). \tag{68}$$

On the other hand, denoting by $r^w(t,x,y)$, $x,y \in \partial\Omega$ the Schwartz kernel of the operator $R^w$ and using the fact that $r^w(\cdot,x,y) = r^w(\cdot,y,x)$ we see that

$$\int_{\partial\Omega} \partial_\nu u^{t^{\alpha+1}F}(x,T)\, H(x)\, dS_x = \int_{\partial\Omega} \partial_\nu u^{t^{\alpha+1}H}(x,T)\, F(x)\, dS_x.$$

Therefore, equation (68) shows that

$$\int_{\partial\Omega} \partial_\nu u^{t^{\alpha+1}F}(x,T)\, H(x)\, dS_x =$$
$$= T^{-\alpha}\left(\Pi^T[t^{\alpha+1}F, t^\alpha H] + \Pi^T[t^\alpha F, t^{\alpha+1}H]\right). \tag{69}$$

Equation (44) follows from (66) and (69).

**Acknowledgements** This work was partly supported by EPSRC (UK) grants GR/M14463 and 36595 and Finnish Academy projects 42013 and 172434. The second and third author are grateful for the hospitality of the Mathematical Sciences Research Institute in Berkeley and support of NSF grant DMS-9810361 (USA). We are grateful to all these organizations.